\documentclass[a4paper,11pt]{amsart}

\usepackage[foot]{amsaddr}

\usepackage{a4wide}
\usepackage{pgfplots}

\usepackage{stmaryrd}

\newcommand{\grad}{\boldsymbol \nabla}

\newcommand{\curl}{\grad \times}

\newcommand{\xx}{\boldsymbol x}

\newcommand{\TT}{\mathcal T}

\newcommand{\FF}{\mathcal F}

\newcommand{\eq}{:=}

\newcommand{\jump}[1]{\left \llbracket #1 \right \rrbracket}
\newcommand{\mean}[1]{\left \{\!\left \{ #1 \right \}\!\right \}}

\newcommand{\ccurl}{\operatorname{curl}}

\newcommand{\JJJ}{\boldsymbol J}
\newcommand{\GGG}{\boldsymbol G}
\newcommand{\EEE}{\boldsymbol E}
\newcommand{\HHH}{\boldsymbol H}

\newcommand{\VVV}{\boldsymbol V}

\newcommand{\MMM}{\boldsymbol M}

\newcommand{\nn}{\boldsymbol n}
\newcommand{\bd}{\boldsymbol d}
\newcommand{\bp}{\boldsymbol p}

\newcommand{\GP}{\Gamma_{\rm P}}
\newcommand{\GA}{\Gamma_{\rm A}}

\newcommand{\pt}{\partial_t}

\newcommand{\eps}{\varepsilon}

\newcommand{\FFP}{\FF^{\rm P}}
\newcommand{\FFA}{\FF^{\rm A}}
\newcommand{\FFI}{\FF^{\rm int}}

\newcommand{\vv}{\boldsymbol v}
\newcommand{\ww}{\boldsymbol w}
\newcommand{\uu}{\boldsymbol u}

\newcommand{\PP}{\mathcal P}
\newcommand{\PPP}{\boldsymbol{\PP}}

\newcommand{\zero}{\boldsymbol 0}

\title%
[A postprocessing for a DG discretization of Maxwell's equations]%
{A postprocessing technique for a discontinuous Galerkin discretization
of time-dependent Maxwell's equations}

\author{G. Nehmetallah$^{\star,\dagger}$}
\author{T. Chaumont-Frelet$^{\star,\dagger}$}
\author{S. Descombes$^{\dagger,\star}$}
\author{S. Lanteri$^{\star,\dagger}$}

\address{\vspace{-.5cm}}
\address{\noindent \tiny \textup{$^\star$Inria, 2004 Route des Lucioles, 06902 Valbonne, France}}
\address{\noindent \tiny \textup{$^\dagger$Laboratoire J.A. Dieudonn\'e, Parc Valrose, 28 Avenue Valrose, 06108 Nice Cedex 02, 06000 Nice, France}}

\begin{document}

\maketitle
\thispagestyle{empty}

\begin{abstract}
We  present  a  novel  postprocessing technique  for  a  discontinuous
Galerkin  (DG) discretization  of  time-dependent Maxwell's  equations
that we couple with an  explicit Runge-Kutta time-marching scheme. The
postprocessed electromagnetic  field converges  one order  faster than
the  unprocessed  solution  in  the  $H(\ccurl)$-norm.   The  proposed
approach is local, in the sense that the enhanced solution is computed
independently in each cell of the computational mesh, and at each time
step  of  interest.   As  a  result, it  is  inexpensive  to  compute,
especially if the  region of interest is localized, either  in time or
space. The  key ideas behind  this postprocessing technique  stem from
hybridizable   discontinuous  Galerkin   (HDG)   methods,  which   are
equivalent  to  the  analyzed  DG   scheme  for  specific  choices  of
penalization parameters. We present several numerical experiments that
highlight  the   superconvergence  properties  of   the  postprocessed
electromagnetic field approximation.

\vspace{.5cm}
\end{abstract}

% \newpage

\section{Introduction}
\label{intro}

Maxwell's  equations  are the  most  general  model of  electrodynamic
theory \cite{griffiths_1999a}.   As a result,  they are employed  in a
variety  of applications,  ranging from  telecommunication engineering
\cite{russer_2006a} to nanophotonics  \cite{gaponenko_2010a}, to study
the propagation  of an electromagnetic field and  its interaction with
structures and matter.

Nowadays,  numerical schemes  are routinely  employed to  simulate the
propagation   of  electromagnetic   waves  by   computing  approximate
solutions              to             Maxwell's              equations
\cite{bondeson_rylander_ingelstrom_2013a}.  While  several approaches,
such as finite difference  methods \cite{yee_1966a}, are available, we
focus      here      on       discontinuous      Galerkin      methods
\cite{fezoui_lanteri_lohrengel_piperno_2005a,hesthaven_warburton:2002,viquerat_2015a},
which have  recently received a lot  of attention, due to  their great
flexibility and ability to handle complex geometries.

Even if currently available computational  power allows for useful and
realistic   simulations,  modeling   accurately  the   propagation  of
electromagnetic fields  in complex geometries is  still a challenging 
and very costly  task. As a result, numerical schemes  are expected to
be accurate and robust, but also  very efficient and adapted to modern
computer architectures.

In the  context of  finite element methods,  postprocessing techniques
are an attractive  way to improve the accuracy of  an already computed
discrete approximation.  In many  cases, these techniques can increase
the order  of convergence of  the method at  a very moderate  cost. In
addition, they  often have  a ``local'' nature,  which allows  for the
design  of  embarrassingly  parallel implementations.   As  a  result,
postprocessing  techniques  and   superconvergence  have  attracted  a
considerable       attention      in       the      past       decades
\cite{Andreev-etal:1988,Babuska-etal:1994,Babuska-etal:1996,Goodsell-etal:1989}.

In  this work,  we elaborate  a novel postprocessing  technique for
time-dependent Maxwell's  equations.  Following \cite{viquerat_2015a},
Maxwell's equations  are discretized with a  first-order discontinuous
Galerkin method coupled with  an explicit Runge-Kutta time-integration
scheme \cite{carpenter-kennedy:1994}. This postprocessing improves the
convergence rate in the $H(\ccurl)$-norm by one order. As with similar
postprocessing techniques  devised in the past,  our proposed approach
is  local,  in  the  sense  that the  enhanced  solution  is  computed
independently in each cell of the computational mesh, and at each time
step of interest.  This is a key property as (a) it enables the design
of highly  parallel numerical  algorithms, and  (b) when  the targeted
application only  requires the knowledge of  the electromagnetic field
in a limited  region of space and/or time, the  amount of computations
is greatly reduced.   Our postprocessing technique is  inspired by two
recent   works,  namely,   a  postprocessing   for  an   explicit  HDG
discretization     of     the     2D    acoustic     wave     equation
\cite{stanglmeier-etal:2016},   and  a   postprocessing   for  a   HDG
discretization   of   the   3D   time-harmonic   Maxwell's   equations
\cite{abgrall_shu:2016}.

We  do  not  carry  out  the mathematical  analysis  of  the  proposed
postprocessing  but   instead,  we  present  a   number  of  numerical
experiments highlighting its  main features. As a result,  our work is
organized as  follows:  in  Section \ref{sec:Settings}, we  recall the
settings   and  key   notations   related   to  Maxwell's   equations,
discontinuous Galerkin methods, and  Runge-Kutta schemes.  We describe
our  postprocessing in  Section  \ref{sec:postprocessing}, and  Section
\ref{sec:Num_exp}  presents numerical  illustrations of  the resulting
methodology.

\section{Settings}
\label{sec:Settings}

\subsection{Maxwell's equations}
\label{subsec:Maxwell}

\begin{subequations}
\label{eq_maxwell_strong}
  
We consider Maxwell's  equations set in a  Lipschitz polyhedral domain
$\Omega   \subset    \mathbb   R^3$    and   in   a    time   interval
$(0,T)$. Specifically,  given $\JJJ:  (0,T) \times \Omega  \to \mathbb
R^3$, the  electromagnetic field $\EEE,\HHH: (0,T)  \times \Omega \to
\mathbb R^3$ satisfies
\begin{equation}
\label{eq_maxwell_volume}
\left \{
\begin{array}{rcl}
\eps \pt \EEE - \curl \HHH & = & \JJJ, \\
\mu \pt \HHH + \curl \EEE  & = & \boldsymbol 0,
\end{array}
\right .
\end{equation}
in $(0,T)  \times \Omega$, where  the functions $\eps,\mu:  \Omega \to
\mathbb R$  respectively represent  the electric permittivity  and the
magnetic  permeability  of the  materials  contained  in $\Omega$.  We
assume that $0 < c \leq \eps,\mu  \leq C$ a.e. in $\Omega$ for fixed
constants $c$ and $C$.

The boundary of $\Omega$ is split into two subdomains $\GA$ and $\GP$,
and we prescribe the boundary conditions
\begin{equation}
\label{eq_maxwell_boundary}
\left \{
\begin{array}{rcll}
\EEE \times \nn_\Omega + \sqrt{\dfrac{\mu}{\eps}} (\HHH \times \nn_\Omega) \times \nn_\Omega 
& = & \GGG & \text{ on } (0,T) \times \GA, \\
\EEE \times \nn_\Omega & = & \boldsymbol 0 & \text{ on } (0,T) \times \GP,
\end{array}
\right .
\end{equation}
where $\nn_\Omega$ denotes the unit vector normal to $\partial \Omega$
pointing outward $\Omega$ and $\GGG: (0,T) \times \GA \to \mathbb R^3$
is a  tangential load term  (i.e.  $\GGG  \cdot \nn_\Omega =  0$). The
first  relation   of  \eqref{eq_maxwell_boundary}  is   a  first-order
absorbing boundary condition (ABC) known  as the Silver-Muller ABC. It
is the  simplest form of  ABC for  Maxwell's equations, and  one could
alternatively consider  higher order ABCs  \cite{hagstrom_lau:2007} or
perfectly matched layers \cite{viquerat_2015a}. The second equation in
\eqref{eq_maxwell_boundary}  models   the  boundary  of   a  perfectly
conducting  material.   Finally,  initial conditions  are  imposed  in
$\Omega$
\begin{equation}
\label{eq_maxwell_initial}
\left \{
\begin{array}{rcl}
\EEE|_{t=0} &=& \EEE_0,
\\
\HHH|_{t=0} &=& \HHH_0,
\end{array}
\right .
\end{equation}
where $\EEE_0,\HHH_0: \Omega \to \mathbb R^3$ are given functions.

Classically \cite{assous_ciarlet_labrunie_2018a}, under the assumption
that  the  data $\mu$,  $\eps$,  $\JJJ$,  $\GGG$, $\EEE_0$  and
$\HHH_0$  are  sufficiently smooth,  there  exists  a unique  pair  of
solution $(\EEE,\HHH)$ to \eqref{eq_maxwell_strong}.
\end{subequations}

We finally mention that in many applications, $\GGG$ is defined
in order to inject an ``incident'' field in the domain. In this
case, we have
\begin{equation}
\label{eq_GGG}
\GGG \eq \EEE^{\rm inc} \times \nn_\Omega +
\sqrt{\dfrac{\mu}{\eps}} (\HHH^{\rm inc} \times \nn_\Omega) \times \nn_\Omega, 
\end{equation}
where $(\EEE^{\rm inc},\HHH^{\rm inc})$ is  a solution to Maxwell's
equations in free space. An important example
that we will  consider in Section \ref{sec:Num_exp} is  the case where
the incident field is a plane wave.

\subsection{Mesh and notations}
\label{subsec:notations}

The domain  $\Omega$ is  partitioned into a  mesh $\TT_h$.   We assume
that  $\TT_h$  consists  of  straight tetrahedral  elements  $K$,  but
hexahedral and/or  curved elements  could be  considered as  well.  We
assume that $\varepsilon$ and  $\mu$ take constant  values $\varepsilon_K$
and $\mu_K$ in each element $K \in \TT_h$.

For the sake  of simplicity, we restrict our attention  to meshes that
are conforming in the  sense of \cite{ciarlet:2002}. Specifically, the
intersection  $\overline{K}_-  \cap  \overline{K}_+$ of  two  distinct
elements $K_\pm \in  \TT_h$ is either a  full face, a full  edge, or a
single  vertex of both $K_-$ and $K_+$. In  particular, hanging  nodes are  not
covered by the  present analysis. This is not  an intrinsic limitation
of the method, but this  assumption greatly simplifies the forthcoming
presentation.

We  denote by  $\FF_h$  the  faces of  the  partition. Recalling  that
$\TT_h$  is  conforming,  each  face  $F  \in  \FF_h$  is  either  the
intersection $\partial K_-  \cap \partial K_-$ of  two elements $K_\pm
\in  \TT_h$, or  is contained  in  the intersection  $\partial K  \cap
\partial \Omega$ of  a single element $K \in \TT_h$  with the boundary
of the domain.
%We also assume that for each face $F \in \FF_h\cap
%\partial \Omega$ either
%$F  \subset  \GP$ or  $F  \subset  \GA$.
   We respectively  denote  by
$\FFI_h$, $\FFP_h$ and  $\FFA_h$ the set internal faces,  and the sets
of faces belonging to $\GP$ and $\GA$.

We associate with each face $F  \in \FF_h$ a unit normal $\nn_F$, with
the  convention  that $\nn_F  =  \nn_\Omega$  if  $F \in  \FFP_h  \cup
\FFA_h$.   If  $F  \in  \FFI_h$,  the orientation  of  the  normal  is
arbitrary, but fixed.  If $\vv: \Omega  \to \mathbb R^3$ is a function
admitting  well-defined  traces  on   $F  \in  \FF_h$,  the  notations
$\jump{\vv}_F$ and $\mean{\vv}_F$ denote the ``jump'' and the ``mean''
of $\vv$ on $F$.  If $F \in  \FFI_h$ with $F = \partial K_- \cap K_+$,
these quantities are defined by
\begin{equation*}
\jump{\vv}_{F} \eq \vv_+|_F (\nn_+ \cdot \nn_F) + \vv_-|_F (\nn_- \cdot \nn_F),
\qquad
\mean{\vv}_{F} \eq \frac{1}{2} \left (\vv_+|_F + \vv_-|_F \right ), 
\end{equation*}
where  $\vv_\pm  \eq  \vv|_{K_\pm}$  and $\nn_\pm$  denotes  the  unit
outward normal to $K_\pm$, while we simply set
\begin{equation*}
\jump{\vv}_{F} \eq \mean{\vv}_{F} \eq \vv|_F, 
\end{equation*}
if $F \in \FFP_h \cup \FFA_h$.

In  the  remaining  of this  work,  $k $  is a  fixed non-negative  integer
representing a  polynomial degree.  For  every element $K  \in \TT_h$,
$\PP_k(K)$ denotes  the set  of polynomials defined  on $K$  of degree
less than or  equal $k$, and $\PPP_k(K) \eq  (\PP_k(K))^3$ denotes the
space  of vector-valued  functions having  polynomial components.   We
finally employ the notation
\begin{equation*}
\PPP_k(\TT_h) \eq \left \{
\vv: \Omega \to \mathbb R^3 \; | \; \vv|_K \in \PPP_k(K) \; \forall K \in \TT_h
\right \}, 
\end{equation*}
for the  space of piecewise  polynomial functions. We also  employ the
notation  $\PPP_k^{\rm  t}(F)$ for the set  of  vector-valued  polynomial
functions defined  on $F$  that are  tangential to  $F$.  $\PPP_k^{\rm
  t}(\FF_h)$ is then  the set of tangential polynomial  defined on the
skeleton of the mesh that are piecewise in $\PPP_k^{\rm t}(F)$.

% If  $\MMM_h^{\rm  t}  \in   \PPP_k^{\rm  t}(\FF_h)$,  and  $\vv_h  \in
% \PPP_k(\TT_h)$, we introduce the notation
% \begin{equation*}
% \langle \MMM_h^{\rm t},\nn \times \vv \rangle_{\FF_h}
% \eq
% \sum_{F \in \FF_h} \int_F \MMM_h^{\rm t} \cdot (\nn_F \times \vv_h).
% \end{equation*}

\subsection{The discontinuous Galerkin scheme}
\label{subsec:DG}

We seek the discrete fields  as piecewise polynomial functions, namely
$\EEE_h,\HHH_h \in  \PPP_k(\TT_h)$. Following \cite{arnold-etal:2002},
the first  step is to  multiply \eqref{eq_maxwell_volume} by  two test
functions $\vv$ and $\ww$, and integrate by parts over each element $K
\in \TT_h$.  We obtain
\begin{equation}
\label{eq_dg_semi_discrete}
\left \{
\begin{array}{rcl}
(\eps \pt \EEE_h,\vv)_{\TT_h}
-
(\HHH_h,\curl \vv)_{\TT_h}
+
\langle \widehat \HHH_h^{\rm t}, \jump{\vv} \times \nn \rangle_{\FF_h}
& = &
(\JJJ,\vv), \\
(\mu \pt \HHH_h,\ww)_{\TT_h}
+
(\EEE_h,\curl \ww)_{\TT_h}
-
\langle \widehat \EEE_h^{\rm t}, \jump{\ww} \times \nn \rangle_{\FF_h}
& = &
0,
\end{array}
\right .
\end{equation}
where $\widehat \EEE_h^{\rm t},  \widehat \HHH_h^{\rm t} \in \PPP^{\rm t}_k(\FF_h)$
are face-based  tangential  fields called  ``numerical fluxes'', and
\begin{equation*}
\langle \widehat \MMM_h^{\rm t}, \jump{\uu} \times \nn \rangle_{\FF_h}
\eq
\sum_{F \in \FF_h}
\int_F 
\widehat \MMM_h^{\rm t} \cdot \left (\jump{\uu}_F \times \nn_F \right ),
\end{equation*}
for  $\MMM_h^{\rm  t}  \in  \PPP_k^{\rm  t}(\FF_h)$,  and  $\uu_h  \in
\PPP_k(\TT_h)$.   We make  use of  numerical fluxes  in the  spirit of
local    DG   methods    that    were    originally   introduced    in
\cite{castillo-etal:2000}  for   scalar elliptic  equations,  and   later  in
\cite{hesthaven_warburton:2002}  for Maxwell's  equations.  We  follow
\cite{viquerat_2015a} to  define our numerical  fluxes.  Specifically,
we set $Z_K  \eq \sqrt{\mu_K/\eps_K}$ and $Y_K \eq 1/Z_K$  for each $K
\in \TT_h$, and we select
\begin{align*}
\widehat \EEE_h^{\rm t}|_F
& \eq
\frac{1}{\mean{Y}}
\left (
\mean{Y \EEE_h}_F^{\rm t} + \dfrac{1}{2} \jump{\HHH_h}_F \times \nn
\right ),
\\
\widehat \HHH_h^{\rm t}|_F
& \eq
\frac{1}{\mean{Z}}
\left (
\mean{Z \HHH_h}_F^{\rm t} - \dfrac{1}{2} \jump{\EEE_h}_F \times \nn
\right ),
\end{align*}
for all $F = \partial K_- \cap \partial K_+ \in \FFI_h$, and 
\begin{equation*}
\widehat \EEE_h^{\rm t}|_F
\eq
\boldsymbol 0
\qquad
\widehat \HHH_h^{\rm t}|_F
\eq
-Y  \EEE_h \times \nn + \HHH_h^{\rm t},
\end{equation*}
if $F = \partial K \cap \GP \in \FFP_h$, and
\begin{align*}
\widehat \EEE_h^{\rm t}|_F
& \eq
\frac{1}{2}
\left (\EEE_h^{\rm t}
+
Z \HHH_h \times \nn  + \GGG \times \nn
\right ),
\qquad
\\
\widehat \HHH_h^{\rm t}|_F
& \eq
\frac{Y}{2}
\left (
Z \HHH_h^{\rm t}  -  \EEE_h \times \nn - \GGG
\right ),
\end{align*}
when $F = \partial K \cap \GA \in \FFA_h$.

\subsection{Time discretization}
\label{subsec:time_disc}

We  can  rewrite  problem \eqref{eq_dg_semi_discrete}  obtained  after
space discretization as
\begin{equation}
\label{eq_ODE_mass}
M \dot U_h(t) + KU_h(t) = B(t), \quad U_h(0) = U_{h,0}
\end{equation}
where  for  each $t  \in  [0,T]$,  the  vector $U_h(t)$  contains  the
coefficients defining  $\EEE_h(t)$ and $\HHH_h(t)$ in  the nodal basis
of  $\PPP_k(\TT_h)$, $M$  and $K$  are  the usual  mass and  stiffness
matrices associated with \eqref{eq_dg_semi_discrete}, and $U_{h,0}$ is
the  interpolation of  the  initial conditions  in the  discretization
space.

Classically, the  key asset of DG  schemes is that the  mass matrix is
block-diagonal, and hence, easy to invert. Thus, we may safely rewrite
\eqref{eq_ODE_mass} as
\begin{equation}
\label{eq_ODE}
\dot U_h(t) = - GU_h(t) + F(t), \quad U_h(0) = U_{h,0},
\end{equation}
where $G \eq M^{-1}  K$ and $F(t) \eq M^{-1} B(t)$.  At this point, we
recognize  in   \eqref{eq_ODE}  a  system  of   ordinary  differential
equations that can be discretized with a time marching scheme.

Here, we focus on a low storage Runge-Kutta scheme, usually denoted by
LSRK($5$,$4$),  presented  in
\cite{carpenter-kennedy:1994}.  After  fixing a time-step  $\Delta t$,
we iteratively  construct approximations  $U_h^n$ of  $U_h(t_n)$, $t_n
\eq n\Delta t$.  Specifically, we let $U^0_h \eq U_{h,0}$,  and for $n
\geq 0$,  $U_h^{n+1}$ is  deduced from  $U_h^n$ through  the following
algorithm
\begin{equation*}
\left\{
\begin{array}{l}
V_h^1 = U_h^{n} \\ \\ 
\left .   
\begin{array}{lcl}  
V_h^{2}
& = &
a_{k}V_h^2
+
\Delta t
\left (
G V_h^1 + F(t_n + c_k \Delta T)
\right )
\\ [0.25cm]
V_h^{1}
& = &
V_h^1 + b_{k}V_h^2
\end{array}
\right\}  
~\hbox{for}~ \,\, k=1,\cdots, 5 \\ \\
U_h^{n+1} = V_h^{1},
\end{array}
\right . 
\end{equation*}
where the coefficients  $a_k$, $b_k$ and $c_k$ are  described in Table
\ref{tab:rk_coeff}. Then, $\EEE_{h,n}$ and $\HHH_{h,n}$ are the element
of $\PPP_k(\TT_h)$ expended on the nodal basis with the coefficients
stored in $U^n_h$.

The  above scheme  is of  particular  interest as  it is  fourth-order
accurate with respect  to the time step $\Delta t$  while being memory
efficient. Indeed,  it only  requires the  storage of  two coefficient
vectors in memory.

Classically, as this time integration scheme is explicit, it is
stable under a CFL condition  linking together the mesh size $h$
and the  selected time step  $\Delta t$. Specifically, given  a mesh
$\TT_h$, we fix the time step by
\begin{equation}
\label{eq_cfl}
\Delta t \eq \alpha_k \min_{K \in \TT_h} \frac{1}{c_K} \frac{V_K}{A_K}
\end{equation}
where, $c_K \eq 1/\sqrt{\varepsilon_K\mu_K}$ is  the wave speed in the
element $K$, and  $V_K$ and $A_K$ are respectively the  volume and the
area  of $K$.  The constant  $\alpha_k$ is  selected according  to the
polynomial degree  $k$. Here, we use the values listed  in Table \ref{table_alpha_cfl},
that we obtained after testing the scheme on simple test-cases.
\begin{table}
\begin{tabular}{llllll}
\hline\noalign{\smallskip}
Coeff  & Value & Coeff & Value  & Coeff & Value \\
\noalign{\smallskip}\hline\noalign{\smallskip}
$a_{1}$ & 0 & $b_{1}$ & $\dfrac{1432997174477}{9575080441755}$ & $c_{1}$ & 0
\\
\noalign{\smallskip}\hline\noalign{\smallskip}	                                                        
$a_{2}$ & $-\dfrac{567301805773}{1357537059087}$ & $b_{2}$ & $\dfrac{5161836677717}{1361206829357}$ & $c_{2}$ & $\dfrac{1432997174477}{9575080441755}$ 
\\
\noalign{\smallskip}\hline\noalign{\smallskip}	                                                          
$a_{3}$ & $-\dfrac{2404267990393}{2016746695238}$ & $b_{3}$ & $\dfrac{1720146321549}{2090206949498}$
& $c_{3}$ & $\dfrac{2526269341429}{6820363962896}$
\\
\noalign{\smallskip}\hline\noalign{\smallskip}	                                                        
$a_{4}$ & $-\dfrac{3550918686646}{2091501179385}$ & $b_{4}$ & $\dfrac{3134564353537}{4481467310338}$ & $c_{4}$ & $\dfrac{2006345519317}{3224310063776}$ \\
\noalign{\smallskip}\hline\noalign{\smallskip}	                                                        
$a_{5}$ & $-\dfrac{1275806237668}{842570457699}$ & $b_{5}$ & $\dfrac{2277821191437}{14882151754819}$
                     & $c_{5}$ & $\dfrac{2802321613138}{2924317926251}$ \\
\noalign{\smallskip}\hline
\end{tabular}
\caption{Values of the coefficients of the LSRK(5,4) scheme.}
\label{tab:rk_coeff}  
\end{table}

\begin{table}
\centering
\begin{tabular}{|c|cccc|}
\hline
$k$ & 1 & 2 & 3 & 4
\\
\hline
$\alpha_k$ & 0.70 & 0.46 & 0.30 & 0.21
\\
\hline
\end{tabular}
\caption{Values of $\alpha_k$ in CFL condition \eqref{eq_cfl}.}
\label{table_alpha_cfl}
\end{table}

Finally, to  ease the discussions  in numerical experiments  below, we
denote by $N$ the number of time steps performed in each simulations.

\section{A novel postprocessing}
\label{sec:postprocessing}

As  discussed above,  $\EEE_{h,n}$ and  $\HHH_{h,n}$ are  respectively
meant to approximate $\EEE(t_n)$ and  $\HHH(t_n)$. The purpose of this
section is to introduce postprocessed solutions $\EEE_{h,n}^\star$ and
$\HHH_{h,n}^\star$   that  are   more   accurate  representations   of
$\EEE(t_n)$ and  $\HHH(t_n)$. This  postprocessing is purely  local in
time,  in the  sense that  the computation  of $\EEE_{h,n}^\star$  and
$\HHH_{h,n}^\star$ only involves $\EEE_{h,n}$  and $\HHH_{h,n}$. It is
also local  in space as the  computation are local to  each element $K
\in          \TT_h$.          Actually,          $\EEE_{h,n}^\star|_K$
(resp. $\HHH_{h,n}^\star|_K$) only depends on $\EEE_{h,n}|_{\widetilde
  K}$ (resp.   $\HHH_{h,n}|_{\widetilde K}$), where $\widetilde  K$ is
the union of  all elements $K' \in \TT_h$ sharing  (at least) one face
with $K$.

Our  approach closely  follows previous  works. Specifically,  similar
postprocessing  strategies have  been  derived  for the  time-harmonic
Maxwell's equations \cite{abgrall_shu:2016}, as well as time-dependent
acoustic  wave  equation  \cite{stanglmeier-etal:2016}.   These  works
develop in  the context  of hybridizable discontinuous  Galerkin (HDG)
methods,  but   can  be  easily   applied  to  the  DG   scheme  under
consideration, as we depict hereafter.

%\subsection{Definition of the post-processed solution}
\label{subsec:Def_post_proc}

Our postprocessing hinges on  element-wise finite element saddle-point
problems.  For each element $K \in  \TT_h$, there exists a unique pair
$(\EEE_{h,n}^\star,p) \in \PPP_{k+1}(K)  \times \PP_{k+2}(K) / \mathbb
R$ such that
\begin{equation*}
\left \{
\begin{array}{rcll}
(\curl \EEE_{h,n}^\star,\curl \ww)_K + (\grad p,\ww)_K
& = &
(\curl \EEE_{h,n},\curl \ww)_K
\\
& + &
\langle \EEE_{h,n}^{\rm t} -
\widehat \EEE_{h,n}^{\rm t}, \nn \times \curl \ww \rangle_{\partial K},
\\
(\EEE_{h,n}^\star,\grad v)_K
& = &
(\EEE_{h,n},\grad v)_K,
\end{array}
\right .
\end{equation*}
for all $\ww \in \PPP_{k+1}(K)$ and  $v \in \PP_{k+2}(K) / \mathbb R$.
Similarly,  for  the  magnetic  field,  there  exists  a  unique  pair
$(\HHH_{h,n}^\star,q)  \in \PPP_{k+1}(K)  \times \PP_{k+2}(K) /
  \mathbb R$ such that
\begin{equation*}
\left \{
\begin{array}{rcll}
(\curl \HHH_{h,n}^\star,\curl \ww)_K + (\grad q,\ww)_K
& = &
(\curl \HHH_{h,n},\curl \ww)_K
\\
& + &
\langle \HHH_{h,n}^{\rm t}-\widehat \HHH_{h,n}^{\rm t}, \nn \times \curl \ww \rangle_{\partial K},
\\
(\HHH_{h,n}^\star,\grad v)_K
& = &
(\HHH_{h,n},\grad v)_K,
\end{array}
\right .
\end{equation*}
for all $\ww \in \PPP_{k+1}(K)$ and  $v \in \PP_{k+2}(K) / \mathbb R$.
$\EEE_{h,n}^\star$ and  $\HHH_{h,n}^\star$ are then  our postprocessed
approximations to $\EEE(t_n)$ and $\HHH(t_n)$.

The left-hand  sides of the  above definition lead to  solve symmetric
linear  systems  of  small  size.  In  addition,  observing  that  the
left-hand  side  is  actually  the same  for  the  two  postprocessing
schemes, we deduce that only  one matrix factorization is required per
element.

The right-hand  sides further show  that for  each $K \in  \TT_h$, the
postprocessed    field   $\EEE_{h,n}^\star|_K$    only   depends    on
$\EEE_{h,n}|_K$ and  the value  at the flux  $\widehat \EEE_{h,n}^{\rm
  t}|_F$  on each  face $F  \in  \FF_K$. In  turn, since  the flux  is
defined  using the  two elements  sharing the  face $F$,  we see  that
$\EEE_{h,n}|_K$ depends on the values taken by $\EEE_{h,n}$ on all the
elements $K'$  sharing at least one  face with $K$. A  similar comment
holds true for $\HHH_{h,n}^\star$.

\section{Numerical experiments}
\label{sec:Num_exp}

\subsection{Standing wave in a cavity}
\label{subsec:CC}

We first consider a model problem given by the propagation of standing
wave in unit cube $\Omega \eq  (0,L)^{3}$, $L \eq 1$ m, with perfectly
conducting  walls   (i.e.   $\GP  \eq   \partial  \Omega$  and   $\GA  \eq
\emptyset$).    Specifically,   we    consider   Maxwell's   equations
\eqref{eq_maxwell_strong}  with  right-hand  sides $\JJJ  \eq  \zero$,
$\GGG \eq \zero$ and initial conditions
\begin{equation*}
\EEE|_{t=0}
\eq
\left (
\begin{array}{r}
-\cos(\pi\xx_1)\sin(\pi\xx_2)\sin(\pi\xx_3)
\\
0 \qquad \qquad \qquad
\\
 \sin(\pi\xx_1)\sin(\pi\xx_2)\cos(\pi\xx_3)
\end{array}
\right ),
\end{equation*}
and $\HHH|_{t=0} \eq \zero$.  $\varepsilon$ and $\mu$ are respectively
set to the vacuum values  $\varepsilon_0 \eq (1/36\pi) \times 10^{-9}$
Fm$^{-1}$ and $\mu_0 \eq 4\pi \times 10^{-7}$ Hm$^{-1}$, and we select
the  simulation  time  $T  \eq  10$ ns.  The  analytical  solution  is
available, and reads
\begin{equation*}
% \left\{
% \begin{array}{lcl}
\EEE(t,\xx) \eq \cos(\omega t)\left (
\begin{array}{r}
-\cos(\pi \xx_1)\sin(\pi\xx_2)\sin(\pi\xx_3)
\\
0 \qquad \qquad \qquad
\\
\sin(\pi \xx_1)\sin(\pi\xx_2)\cos(\pi\xx_3)
\end{array}
\right ),
\end{equation*}
and
\begin{equation*}
\HHH(t,\xx) \eq
\frac{\pi}{\omega}\sin(\omega t)
\left (
\begin{array}{r}
 \sin(\pi\xx_1)\cos(\pi\xx_2)\cos(\pi\xx_3)
\\
2\cos(\pi\xx_1)\sin(\pi\xx_2)\cos(\pi\xx_3)
\\
 \cos(\pi\xx_1)\cos(\pi\xx_2)\sin(\pi\xx_3)
\end{array}
\right ),
\end{equation*}
where  the  angular frequency  is  given  by $\omega  \eq  \sqrt{3}\pi
c_0/L$,  $c_0  \eq  1/\sqrt{\varepsilon_0\mu_0}$ being  the  speed  of
light.

We  consider structured  meshes  $\TT_h$ that  are  obtained by  first
splitting $\Omega$ into $n \times n \times n$ cubes ($n \eq L/h$), and
then splitting each cube into $6$ tetrahedra.

Figures  \ref{fig:Hcurlpperr_E_CC} and  \ref{fig:Hcurlpperr_H_CC} show
the behavior of the error  for the original and postprocessed discrete
solutions  with  respect  to  time  on  a  fixed  mesh  built  from  a
$8\times8\times8$ Cartesian  partition.  The  time step $\Delta  t$ is
selected following CFL condition \eqref{eq_cfl}.  Both  the original
and the postprocessed error exhibit  an oscillatory behavior, which is
typical of  this particular test  case. The postprocessed  solution is
about 10 times more accurate than the original one.

Table \ref{tab:Hcurlpp_order_CC}  presents in more detail  our results
on a  series of meshes and  for different polynomial degrees.   We see
that in  each case, the  curl of the postprocessed  solution converges
with the expected order, namely $k+1$.

\begin{figure}
\centering
%\hspace*{1.5cm}
\begin{tikzpicture}[scale=0.95]
\begin{axis}[ 
legend pos= outer north east,
xlabel={Time (s)}, 
ymajorgrids=true,
xmajorgrids=true,
grid style=dashed,
ymode=log,
ymin=.0001,
ymax=.1,
xmin=0,
xmax=1e-8,
ytick={.1,.01,.001}
]
\addplot[mark=o,draw=blue,thick] table{data/HCURL_Error_E_P2_8_CC};
\addlegendentry{$\|\curl(\EEE(t_n)-\EEE_{h,n})\|_\Omega$}
\addplot[mark=square,draw=red,thick] table{data/HCURLPP_Error_E_P2_8_CC};
\addlegendentry{$\|\curl(\EEE(t_n)-\EEE_{h,n}^\star)\|_\Omega$}
\end{axis}
\end{tikzpicture}
\caption{Standing wave in  a cubic cavity: time evolution  of the error on 
  the electric field.}
\label{fig:Hcurlpperr_E_CC}
\end{figure}
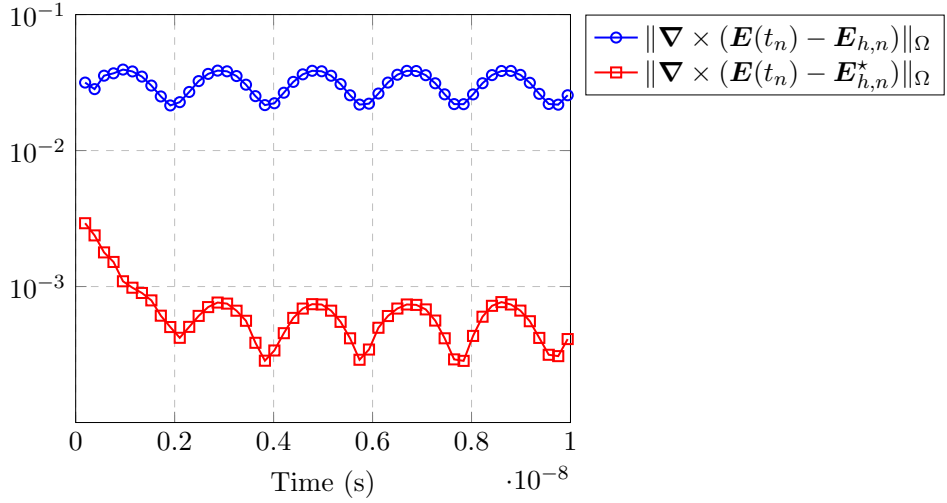
\begin{figure}
\centering
%\hspace*{1.5cm}
\begin{tikzpicture}[scale=0.95]
\begin{axis}[ 
legend pos= outer north east,
xlabel={Time (s)}, 
ymajorgrids=true,
xmajorgrids=true,
grid style=dashed,
ymode=log,
ymin=.0001,
ymax=.1,
xmin=0,
xmax=1e-8,
ytick={.1,.01,.001}
]
\addplot[mark=o,draw=blue,thick] table{data/HCURL_Error_H_P2_8_CC};
\addlegendentry{$\|\curl(\HHH(t_n)-\HHH_{h,n})\|_\Omega$}
\addplot[mark=square,draw=red,thick] table{data/HCURLPP_Error_H_P2_8_CC};
\addlegendentry{$\|\curl(\HHH(t_n)-\HHH_{h,n}^\star))\|_\Omega$}
\end{axis}
\end{tikzpicture}
\caption{Standing wave in  a cubic cavity: time evolution  of the error
  on the magnetic field.}
\label{fig:Hcurlpperr_H_CC}
\end{figure}
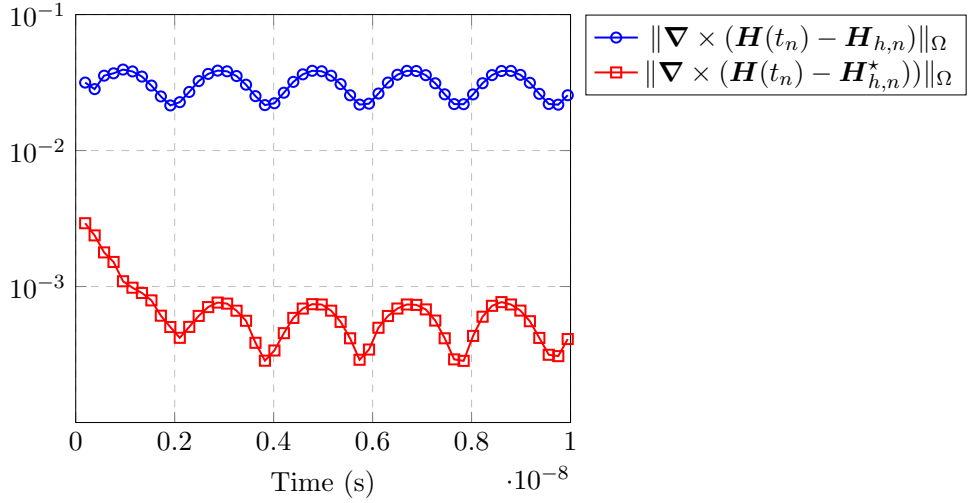
\begin{table}
\label{tab:Hcurlpp_order_CC}   
\begin{tabular}{llll}
\hline\noalign{\smallskip}
& $h$ & $ \|\curl(\EEE(T)-\EEE_{h,N})\|_\Omega$ &
        $\|\curl(\EEE(T)-\EEE_{h,N}^\star)\|_\Omega$ \\
\noalign{\smallskip}\hline\noalign{\smallskip}
         & $1/4$ & 7.99e-01 & 6.37e-01 \\
$\PPP_1$ & $1/6$ & 4.94e-01 \qquad (\textbf{eoc 1.19}) & 2.69e-01 \qquad (\textbf{eoc 2.13}) \\
         & $1/8$ & 3.65e-01 \qquad (\textbf{eoc 1.05}) & 1.45e-01 \qquad (\textbf{eoc 2.15}) \\
\noalign{\smallskip}\hline\noalign{\smallskip}	                                                        
         & $1/4$ & 1.40e-01 & 3.80e-02 \\
$\PPP_2$ & $1/6$ & 6.55e-02 \qquad (\textbf{eoc 1.87}) & 1.04e-02 \qquad (\textbf{eoc 3.20}) \\
         & $1/8$ & 3.75e-02 \qquad (\textbf{eoc 1.94}) & 4.24e-03 \qquad (\textbf{eoc 3.12}) \\
\noalign{\smallskip}\hline\noalign{\smallskip}
         & $1/4$ & 2.05e-02 & 4.32e-03 \\
$\PPP_3$ & $1/6$ & 6.17e-03 \qquad (\textbf{eoc 2.96}) & 9.29e-04 \qquad (\textbf{eoc 3.74}) \\
         & $1/8$ & 2.62e-03 \qquad (\textbf{eoc 2.98}) & 3.09e-04 \qquad (\textbf{eoc 3.83}) \\        
\end{tabular}
\begin{tabular}{llll}
\hline\noalign{\smallskip}
& $h$ & $ \|\curl(\HHH(T)-\HHH_{h,N})\|_\Omega$ &
        $\|\curl(\HHH(T)-\HHH_{h,N}^\star)\|_\Omega$ \\
\noalign{\smallskip}\hline\noalign{\smallskip}
         & $1/4$ & 6.17e-01 & 4.18e-01 \\
$\PPP_1$ & $1/6$ & 3.76e-01 \qquad (\textbf{eoc 1.22}) & 1.80e-01 \qquad (\textbf{eoc 2.08}) \\
         & $1/8$ & 2.70e-01 \qquad (\textbf{eoc 1.15}) & 9.71e-02 \qquad (\textbf{eoc 2.15}) \\
\noalign{\smallskip}\hline\noalign{\smallskip}	                                                        
         & $1/4$ & 9.94e-02 & 2.19e-02 \\
$\PPP_2$ & $1/6$ & 4.68e-02 \qquad (\textbf{eoc 1.86}) & 6.00e-03 \qquad (\textbf{eoc 3.19}) \\
         & $1/8$ & 2.71e-02 \qquad (\textbf{eoc 1.90}) & 2.44e-03 \qquad (\textbf{eoc 3.13}) \\                                                             
\noalign{\smallskip}\hline\noalign{\smallskip}	                                                        
         & $1/4$ & 1.60e-02 & 2.46e-03 \\
$\PPP_3$ & $1/6$ & 4.83e-03 \qquad (\textbf{eoc 2.95}) & 5.39e-04 \qquad (\textbf{eoc 3.74}) \\
         & $1/8$ & 2.06e-03 \qquad (\textbf{eoc 2.96}) & 1.82e-04 \qquad (\textbf{eoc 3.77}) \\
\noalign{\smallskip}\hline
\end{tabular}
\caption{Standing wave  in a  cubic cavity: numerical convergence.}
\end{table}
\subsection{Plane wave in free space}
\label{subsec:PW}

We  now consider  the  propagation  of a  plane  wave  in free  space.
Specifically,       we        consider       Maxwell's       equations
\eqref{eq_maxwell_strong} with  $\Omega \eq  (0,L)^3$, $L
\eq 1$ m, $\GP \eq \emptyset$  and $\GA \eq \partial \Omega$. $\JJJ
\eq \zero$, and $\GGG$ is defined by \eqref{eq_GGG} with
\begin{equation*}
\EEE^{\rm inc}(t,\xx)
\eq
\bp \cos \left (\omega \left (t - \frac{\bd \cdot \xx}{c_0}\right )\right ),
\quad
\HHH^{\rm inc}(t,\xx) \eq \sqrt{\frac{\eps_0}{\mu_0}} \bd \times \EEE^{\rm inc}(t,\xx), 
\end{equation*}
\noindent   where   $\bp \eq (1,0,0)^T$   is   the
polarization, $\bd \eq (0,0,1)^T$ is  the direction of propagation and
$\omega \eq  6 \pi  c_0/L$ is  the angular  frequency.  We  impose the
initial   conditions  \eqref{eq_maxwell_initial}   with  $\EEE_0   \eq
\EEE^{\rm inc}|_{t=0}$ and $\HHH_0  \eq \HHH^{\rm inc}|_{t=0}$.  Then,
since  the medium  under consideration  is homogeneous,  no reflection
and/or diffraction occur, and the  analytical solution is simply $\EEE
=  \EEE^{\rm  inc}$  and  $\HHH  = \HHH^{\rm  inc}$.   We  select  the
simulation  time $T \eq 10$ ns.  As  for the  cubic cavity
test, we consider  structured meshes $\TT_h$, that we  obtain by first
splitting  $\Omega$  into  $n  \times   n  \times  n$  cubes  ($n  \eq
L/h$),    and   then    splitting    each    cube   into    $6$
tetrahedra. As explained above, the  time step is selected using
  \eqref{eq_cfl}.         Figures~\ref{fig:Hcurlpperr_E_PW}       and
\ref{fig:Hcurlpperr_H_PW}  show the  behaviour  of the  error for  the
original and postprocessed discrete solutions  with respect to time on
a fixed mesh based on a $12\times12\times12$ Cartesian partition.  The
postprocessed  solution  is  about  5 times  more  accurate  than  the
original solution.  Table  \ref{tab:Hcurlpp_order_PW} presents in more
detail our results on a series  of meshes and for different polynomial
degrees. We  see that  in each  cases, the  curl of  the postprocessed
solution converges with the expected order, namely $k+1$.
\begin{figure}
\centering
%\hspace*{1.5cm}
\begin{tikzpicture}[scale=0.95]
\begin{axis}[ 
legend pos= outer north east,
xlabel={Time (s)}, 
ymajorgrids=true,
xmajorgrids=true,
grid style=dashed,
ymode=log,
ymin=.1,
ymax=1.5,
xmin=0,
xmax=1e-8,
ytick={1,.1}
]
\addplot[mark=o,draw=blue,thick] table{data/HCURL_Error_E_P2_12_PW};
\addlegendentry{$\|\curl(\EEE(t_n)-\EEE_{h,n})\|_\Omega$}
\addplot[mark=square,draw=red,thick] table{data/HCURLPP_Error_E_P2_12_PW};
\addlegendentry{$\|\curl(\EEE(t_n)-\EEE_{h,n}^\star)\|_\Omega$}
\end{axis}
\end{tikzpicture}
\caption{Plane wave in free space: time evolution of the error on the electric field.} 
\label{fig:Hcurlpperr_E_PW}
\end{figure}
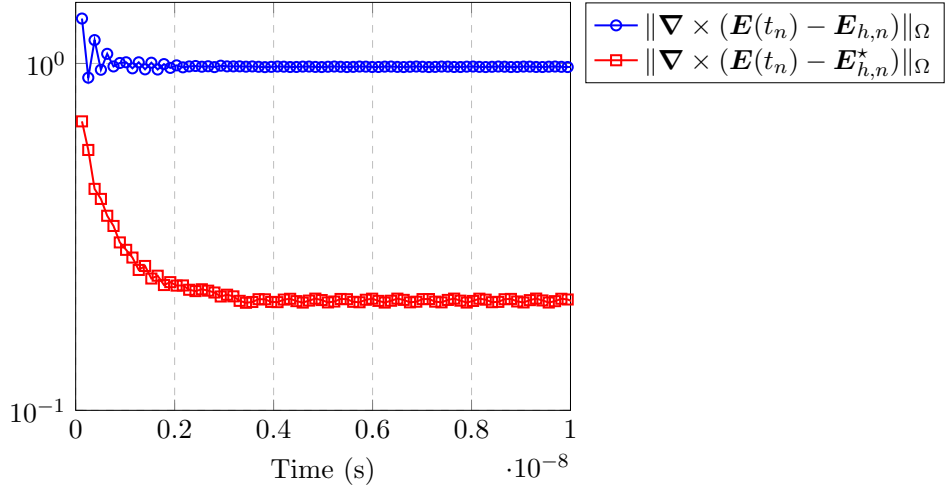
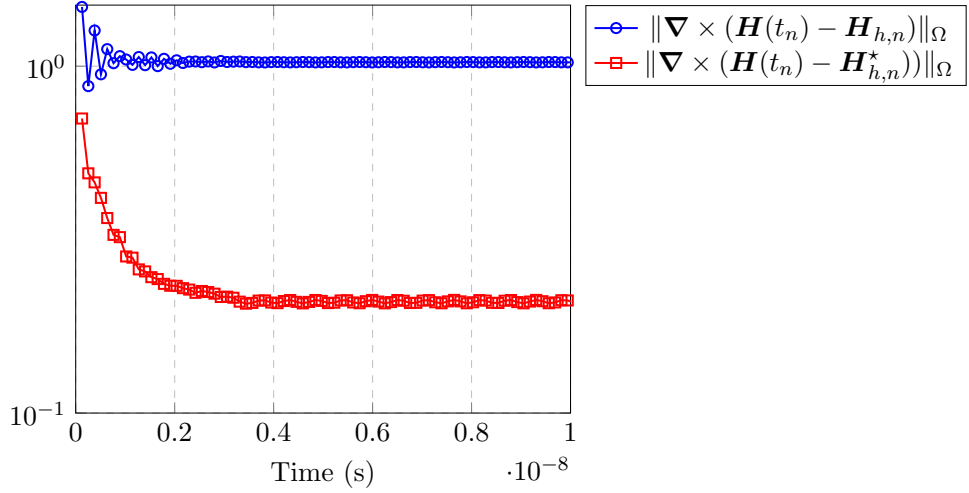
\begin{figure}
\centering
%\hspace*{1.5cm}
\begin{tikzpicture}[scale=0.95]
\begin{axis}[ 
legend pos= outer north east,
xlabel={Time (s)}, 
ymajorgrids=true,
xmajorgrids=true,
grid style=dashed,
ymode=log,
ymin=.1,
ymax=1.5,
xmin=0,
xmax=1e-8,
ytick={1,.1}
]
\addplot[mark=o,draw=blue,thick] table{data/HCURL_Error_H_P2_12_PW};
\addlegendentry{$\|\curl(\HHH(t_n)-\HHH_{h,n})\|_\Omega$}
\addplot[mark=square,draw=red,thick] table{data/HCURLPP_Error_H_P2_12_PW};
\addlegendentry{$\|\curl(\HHH(t_n)-\HHH_{h,n}^\star))\|_\Omega$}
\end{axis}
\end{tikzpicture}
\caption{Plane wave in free space: time evolution of the error on the magnetic field.}
\label{fig:Hcurlpperr_H_PW}
\end{figure}
\begin{table}
\label{tab:Hcurlpp_order_PW}   
\begin{tabular}{llll}
\hline\noalign{\smallskip}
& $h$ & $ \|\curl(\EEE(T)-\EEE_{h,N})\|_\Omega$ &
        $\|\curl(\EEE(T)-\EEE_{h,N}^\star)\|_\Omega$ \\
\noalign{\smallskip}\hline\noalign{\smallskip}
         & $1/8$  & 5.37e-00 & 6.02e-00 \\
$\PPP_1$ & $1/10$ & 4.38e-00 \qquad (\textbf{eoc 0.92}) & 3.99e-00 \qquad (\textbf{eoc 1.84}) \\
         & $1/12$ & 3.75e-00 \qquad (\textbf{eoc 0.86}) & 2.73e-00 \qquad (\textbf{eoc 2.08}) \\
\noalign{\smallskip}\hline\noalign{\smallskip}                  
         & $1/8$  & 1.98e-00 & 7.92e-01 \\
$\PPP_2$ & $1/10$ & 1.36e-00 \qquad (\textbf{eoc 1.70}) & 3.72e-01 \qquad (\textbf{eoc 3.38}) \\
         & $1/12$ & 9.77e-01 \qquad (\textbf{eoc 1.81}) & 2.08e-01 \qquad (\textbf{eoc 3.18}) \\
                                                                
\noalign{\smallskip}\hline\noalign{\smallskip}	                                                        
         & $1/8$  & 4.63e-01 & 1.01e-01 \\
$\PPP_3$ & $1/10$ & 2.44e-01 \qquad (\textbf{eoc 2.88}) & 4.25e-02 \qquad (\textbf{eoc 3.87}) \\
         & $1/12$ & 1.43e-01 \qquad (\textbf{eoc 2.93}) & 2.22e-02 \qquad (\textbf{eoc 3.56}) \\
\end{tabular}
\begin{tabular}{llll}
\hline\noalign{\smallskip}
& $h$ & $ \|\curl(\HHH(T)-\HHH_{h,N})\|_\Omega$ &
        $\|\curl(\HHH(T)-\HHH_{h,N}^\star)\|_\Omega$ \\
\noalign{\smallskip}\hline\noalign{\smallskip}
         & $1/8$  & 5.89e-00 & 6.01e-00 \\
$\PPP_1$ & $1/10$ & 4.68e-00 \qquad (\textbf{eoc 1.03}) & 3.97e-00 \qquad (\textbf{eoc 1.85}) \\
         & $1/12$ & 4.00e-00 \qquad (\textbf{eoc 0.86}) & 2.75e-00 \qquad (\textbf{eoc 2.03}) \\
\noalign{\smallskip}\hline\noalign{\smallskip}	                                                     
         & $1/8$ & 2.16e-00  & 7.60e-01 \\
$\PPP_2$ & $1/10$ & 1.45e-00 \qquad (\textbf{eoc 1.79}) & 3.71e-01 \qquad (\textbf{eoc 3.21}) \\
         & $1/12$ & 1.03e-00 \qquad (\textbf{eoc 1.89}) & 2.11e-01 \qquad (\textbf{eoc 3.10}) \\
\noalign{\smallskip}\hline\noalign{\smallskip}	                                                     
         & $1/8$ & 4.87e-01  & 1.01e-01 \\
$\PPP_3$ & $1/10$ & 2.54e-01 \qquad (\textbf{eoc 2.93}) & 4.32e-02 \qquad (\textbf{eoc 3.79}) \\
         & $1/12$ & 1.48e-01 \qquad (\textbf{eoc 2.96}) & 2.29e-02 \qquad (\textbf{eoc 3.48}) \\
\noalign{\smallskip}\hline
\end{tabular}
\caption{Plane wave in free space: numerical convergence.}
\end{table}

\subsection{Scattering of a plane wave by a dielectric sphere}
\label{subsec:dielectric_sphere}

\begin{figure}
\includegraphics[height=3.0cm,width=1.25cm]{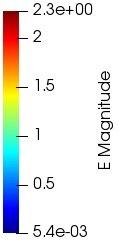}
\includegraphics[height=5.0cm,width=5.25cm]{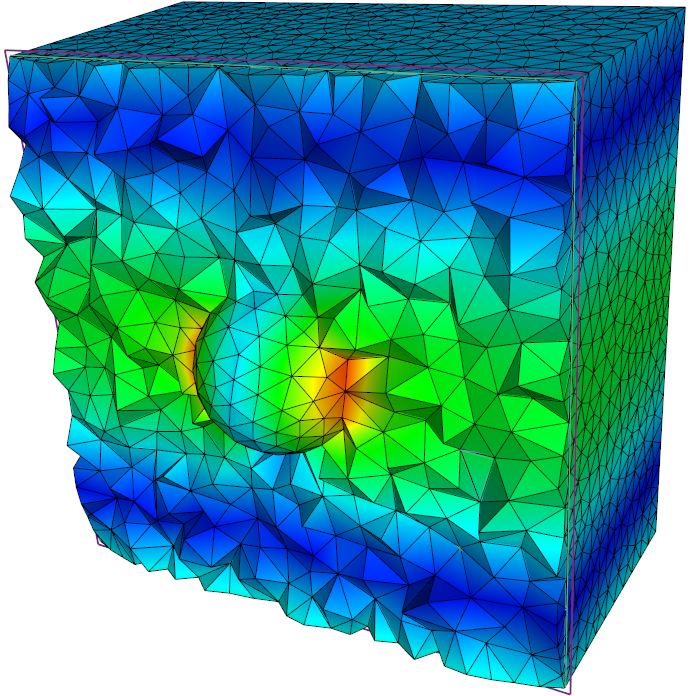}
\caption{Representation of $|\EEE_{\rm r}(T)|$ in the scattering example.}
\label{fig:dielectric_sphere}
\end{figure}

We  now consider  a problem  involving a  dielectric sphere  of radius
0.15~m   with  $\varepsilon=2\varepsilon_0$   and  $\mu=\mu_0$.    The
computational domain  is bounded by  a cube of side  1 m on  which the
Silver-Muller absorbing  condition is applied and  the simulation time
is $T  \eq 3$ ns.   We make use  of an unstructured  tetrahedral mesh,
which consists  of 32,602 elements with  565 elements in the  sphere
and $\Delta  t$ is  chosen via  \eqref{eq_cfl}.  The  right-hand sides
$\JJJ$ and  $\GGG$ are the  same than in Example  \ref{subsec:PW}, and
the initial conditions are taken to be zero.  We select $\PPP_2$
elements,      and     denote      by     $(\EEE_h,\HHH_h)$      and
$(\EEE_h^\star,\HHH_h^\star)$   the   original   and   postprocessed
solutions.    As  the   analytical  solution   to  the   problem  is
unavailable,   we   compute   a   reference   solution   $(\EEE_{\rm
  r},\HHH_{\rm r})$ with $\PPP_4$ elements  on the same mesh and the
time step is defined as $\Delta  t_{\rm r} \eq \Delta t/3$.  $\Delta
t_{\rm  r}$ is  chosen  as an  integral division  of  $\Delta t$  to
facilitate comparisons.  We chose to divide $\Delta t$ by $3$ since,
following Table  \ref{table_alpha_cfl}, it  is the  smallest integer
for which  CFL condition  \eqref{eq_cfl} holds  true.  We  refer the
reader to Figure~\ref{fig:dielectric_sphere} for  a snapshot of the
reference solution.

To assess the impact of the postprocessing, we consider a set of
evaluation points $\boldsymbol A$, and we compute
relative errors
$$
\text{err}(\VVV)^2
=
\frac{\sum_{n=1}^{N} || \curl (\VVV_{\rm r}(t_n,\boldsymbol{A}) - \VVV_{h,n}(\boldsymbol{A}))||^2}%
{\sum_{n=1}^{N} ||\curl (\VVV_{\rm r})(t_n,\boldsymbol{A})||^2}
$$
and
$$
\text{err}^\star(\VVV)^2 =
\frac{\sum_{n=1}^{N} ||\curl (\VVV_{\rm r}(t_n,\boldsymbol{A}) - \VVV_{h,n}^\star(\boldsymbol{A}))||^2}%
{\sum_{n=1}^{N} ||\curl (\VVV_{\rm r})(t_n,\boldsymbol{A})||^2}
$$
with $\VVV \eq \EEE$ or $\HHH$. Table~\ref{tab:err_PP} shows that
our postprocessing approach reduces the error by at least a factor of $2$
for the 9 evaluation points that we have selected.

\begin{table}
\label{tab:err_PP}
\begin{tabular}{lll}
\hline
Point  & Field  &  $\text{err}$ \qquad $\text{err}^\star$ \\ 
\noalign{\smallskip}\hline\noalign{\smallskip}                  
                     & $\EEE$ & \textbf{0.083} \quad \textbf{0.033} \\
$A_1(0,0,0.45)$      & $\HHH$ & \textbf{0.103} \quad \textbf{0.048} \\
\noalign{\smallskip}\hline\noalign{\smallskip}
                     & $\EEE$ & \textbf{0.008} \quad \textbf{0.005} \\
$A_2(0.2,-0.3,0.8)$  & $\HHH$ & \textbf{0.008} \quad \textbf{0.006} \\
\noalign{\smallskip}\hline\noalign{\smallskip}	
                     & $\EEE$ & \textbf{0.019} \quad \textbf{0.005} \\
$A_3(0.2,-0.3,0.2)$  & $\HHH$ & \textbf{0.020} \quad \textbf{0.006} \\
\noalign{\smallskip}\hline\noalign{\smallskip}
                     & $\EEE$ & \textbf{0.015} \quad \textbf{0.004} \\
$A_4(0.2,0.3,0.2)$   & $\HHH$ & \textbf{0.017} \quad \textbf{0.005} \\
\noalign{\smallskip}\hline\noalign{\smallskip}	
                     & $\EEE$ & \textbf{0.019} \quad \textbf{0.007} \\
$A_5(0.2,0.3,0.8)$   & $\HHH$ & \textbf{0.027} \quad \textbf{0.007} \\
\noalign{\smallskip}\hline\noalign{\smallskip}	
                     & $\EEE$ & \textbf{0.015} \quad \textbf{0.008} \\
$A_6(-0.2,-0.3,0.8)$ & $\HHH$ & \textbf{0.014} \quad \textbf{0.008} \\
\noalign{\smallskip}\hline\noalign{\smallskip}
                     & $\EEE$ & \textbf{0.027} \quad \textbf{0.008} \\
$A_7(-0.2,-0.3,0.2)$ & $\HHH$ & \textbf{0.028} \quad \textbf{0.008} \\
\noalign{\smallskip}\hline\noalign{\smallskip}
                     & $\EEE$ & \textbf{0.021} \quad \textbf{0.007} \\
$A_8(-0.2,0.3,0.2)$  & $\HHH$ & \textbf{0.024} \quad \textbf{0.007} \\
\noalign{\smallskip}\hline\noalign{\smallskip}
                     & $\EEE$ & \textbf{0.010} \quad \textbf{0.005} \\
$A_9(-0.2,0.3,0.8)$  & $\HHH$ & \textbf{0.011} \quad \textbf{0.005} \\
\noalign{\smallskip}\hline
\end{tabular}
\caption{Scattering  of a  plane wave  by a  dielectric sphere:  $L^2$
error  between  the  reference  solution and  the  solution  with  a
$\PPP_2$    interpolation   with    and    without   applying    the
postprocessing.}
\end{table}

\section{Conclusion}
\label{sec:conc}

In  this  work we  have  presented  a  postprocessing approach  for  a
discontinuous Galerkin discretization  of the time-dependent Maxwell's
equations in  3D.  This  postprocessing technique is  inexpensive, and
can be computed  independently in each element of  the mesh, and
at every time  step of interest.  It is thus  well adapted to parallel
computer architectures.   Moreover, it is particularly  suited to
applications requiring a higher  accuracy in localized regions, either
in  time or  space. We  have  presented  numerical  examples,  both  with
analytical solution and in  complicated geometries, that indicate that
our  postprocessing  approach improves  the  convergence  rate of  the
discrete solution in the $H(\ccurl)$-norm by one order.  Overall, this
contribution is  to be employed  as an  efficient way of  reducing the
$H(\ccurl)$-norm error of discontinuous Galerkin discretizations.

\bibliographystyle{amsplain}
\bibliography{bibliography.bib}

\end{document}